\documentclass[fleqn,reqno]{amsart}
\usepackage[mathscr]{eucal}
\usepackage{amsfonts}
\usepackage{amsmath}
\usepackage{amsthm}
\usepackage{amssymb}

\newcommand{\C}{\mathbb C}
\newcommand{\Hom}{{\mathrm{Hom}}}
\newcommand{\g}{{\mathfrak{g}}}
\newcommand{\h}{{\mathfrak{h}}}

\newtheorem{theorem}{Theorem}[section]
\newtheorem{theorem/definition}{Theorem/Definition}[section]

\newtheorem{proposition}{Proposition}[section]
\newtheorem{lemma}{Lemma}[section]

\newtheorem{corollary}{Corollary}[section]

\theoremstyle{remark}
\newtheorem{remark}{Remark}[section]
\theoremstyle{definition}
\newtheorem{example}{Example}[section]
\newtheorem{definition}{Definition}[section]

\begin{document}

\title[On the lower central series]{On the lower central series 
of an associative algebra}

\author{Galyna Dobrovolska}
\author{John Kim}
\author{Xiaoguang Ma}

\address{Massachusetts Institute of Technology,
Cambridge, MA 02139, USA}
\email{galyna@mit.edu}

\address{Massachusetts Institute of Technology,
Cambridge, MA 02139, USA}
\email{kimjohn@mit.edu}

\address{Department of Mathematics, Massachusetts Institute of Technology,
Cambridge, MA 02139, USA}
\email{xma@math.mit.edu}


\maketitle

\centerline{\bf (with an appendix by Pavel Etingof)}

\vskip .1in

\begin{abstract}
For an associative algebra $A$, define its lower central series $L_{0}(A)=A$, $L_{i}(A)=[A,L_{i-1}(A)]$, and the corresponding quotients $B_{i}(A)=L_{i}(A)/L_{i+1}(A)$. 
In this paper, we study the structure of $B_{i}(A_{n})$ for a free algebra $A_{n}$. We construct a basis for $B_{2}(A_{n})$ and determine the structure of $B_{3}(A_{2})$ and $B_{4}(A_{2})$. In the appendix, we study the structure of $B_{2}(A)$ for any associative algebra $A$ over $\mathbb{C}$.
\end{abstract}

\section{Introduction}

Let $A$ be an associative algebra. Let us regard it as
a Lie algebra with commutator $[a,b]=ab-ba$. Then one can inductively define 
the lower central series filtration of $A$: 
$L_1(A)=A$, $L_i(A)=[A,L_{i-1}(A)]$, and the corresponding quotients 
$B_i(A)=L_i(A)/L_{i+1}(A)$. It is an interesting problem to understand 
the structure of the spaces $B_i(A)$ for a given algebra $A$. 

The study of $B_i(A)$ was initiated in the paper by B. Feigin and B. Shoikhet
\cite{FS}, who considered the case when $A=A_n$ is the free associative 
algebra in $n$ generators over $\mathbb{C}$. 
Their main results are that 
$B_i(A_n)$ for $i>1$ are representations of the Lie algebra $W_n$
of polynomial vector fields in $n$ variables, 
and that $B_2(A_n)$ is isomorphic, as a
$W_n$-module, to the space of closed (or equivalently, exact) 
polynomial differential forms on $\mathbb{C}^n$ of positive even 
degree. 

The goal of this paper is to continue the study of the structure 
of $B_i(A_n)$, and more generally, 
of $B_i(A)$ for any associative algebra $A$. 
More specifically, in Section 2 we give a new simple proof 
of the result of Feigin and Shoikhet on the structure of
$B_2(A_n)$ for $n=2,3$, by constructing an
explicit basis of this space. In Section 3, we generalize this
basis to the case $n>3$, and use it to determine the structure of the space 
$B_2(A_n^R)$, where $A_n^R$ is the quotient of $A_n$ by the relations   
$x_i^{m_i}=0$, where $m_i$ are positive integers. In Sections 4,5 
we obtain some information about the structure of $B_m(A_2)$ as
a $W_2$-module. In Section 6 we determine the structures of $B_3(A_2)$
and $B_4(A_2)$, thus confirming conjectures from \cite{FS}.
Finally, in the appendix, the structure of $B_2(A)$ is studied for
any associative algebra $A$ over $\Bbb C$.

\section{The structure of $B_{2,2}$ and $B_{3,2}$}
\label{n23}

\subsection{Some notations}
Let $A_n$ be the free algebra over $\Bbb C$ in $n$ generators
$x_{1},\ldots,x_{n}$. Let $m_1,\ldots,m_n$ be positive integers,
and $R$ be the set of relations $x_i^{m_i}=0$, $i=1,\ldots,n$.
Let $A_n^R=A_n/(R)$. From now on, we denote 
$B_{i}(A_{n})$ by $B_{n,i}$, and  $B_{i}(A_{n}^{R})$ by $B_{n,i}^{R}$. 

Let each generator $x_{i}$ have degree $1$.  
We say that $w\in B_{n,k}$ has multidegree $(i_1,\ldots ,i_n)$ if every $x_s$ occurs $i_s$
times in every monomial of $w$. The set of all multidegree ${\mathbf i}=(i_1,\ldots ,i_n)$ elements in $B_{n,k}$ is denoted by $B_{n,k}[{\mathbf i}]$. Note that not all $w\in B_{n,k}$ will have a
multi-degree.  However, monomials and brackets of monomials in $B_{n,k}$ will have a
multi-degree. 
Let $l = i_1+\cdots +i_n$, and  call $l$ the degree of $w$. We denote $B_{n,k}[l]$ to be the set of all 
degree $l$ elements in $B_{n,k}$. 

\subsection{Basis for $B_{2,2}[l]$}
In this section, we find a basis for $B_{2,2}[l]$.

\begin{proposition}
\label{2VarGen}
For $l\geq 2$, the $l-1$ elements $[x_1^i,x_2^{l-i}]$ for $i = 1,\dots,l-1$ constitute
a spanning set for $B_{2,2}[l]$.
\end{proposition}

\begin{proof}
First note that every element of $B_{2,2}[l]$ can be expressed as a linear combination of
the brackets $[a,x_1]$ and $[b,x_2]$, where $a$ and $b$ are monomials with degree
no less than $1$.  To see this, consider an arbitrary bracket of
monomials in $B_{2,2}[l]$.  This bracket may be
written as $[P,q_1q_2\cdots q_n]$, where $n\ge2$ and $q_i$ represents either $x_1$ or $x_2$.

Then we have 
$$
[P,q_1q_2\cdots q_n] = [Pq_1\cdots q_{n-1},q_n]+[q_nPq_1\cdots q_{n-2},q_{n-1}]+\cdots +[q_2\cdots q_nP,q_1].
$$

As every element of $B_{2,2}[l]$ is a linear combination of such brackets of monomials, and
each bracket of monomials is a sum of brackets of the desired form, every element of
$B_{2,2}[l]$ is a linear combination of brackets of the desired form.

Consider $[a,x_1]$.  Write $a = x_1^k a_1$, where $a_1$ begins with $x_2$ or is equal to $1$.
Then we have $[a_1,x_1^{k+1}] = \sum_{j=0}^{k}[x_1^{k-j}a_1x_1^j,x_1]$.  Notice that all the
terms in the summation are equivalent in $B_{2,2}[l]$ because we can cyclically permute either
term of the bracket.  So $[x_1^k a_1,x_1] = \frac{1}{k+1} [a_1,x_1^{k+1}]$.

If $a_1$ is not $1$, then $a_1$ can be written as $x_2^m x_1^n a_2$, where $a_2$ begins with
$x_2$ or is equal to $1$.  So by the same argument as above,
\begin{eqnarray*}
[a_1,x_1^{k+1}] = [x_2^m x_1^n a_2,x_1^{k+1}] =
[x_1^n a_2 x_2^m,x_1^{k+1}] = \frac{1}{n+1} [a_2 x_2^m,x_1^{k+1+n}]
\end{eqnarray*}
Continuing this process will eventually transfer all powers of $x_1$ to the right
side of the bracket, showing that $[a,x_1]$ is a constant multiple of $[x_2^{l-i},x_2^i] =
-[x_1^i,x_2^{l-i}]$ for some $i$ from $1$ to $l-1$.

A similar argument shows that $[b,x_2]$ is a constant multiple of $[x_1^i,x_2^{l-i}]$.
Recalling that every element of $B_{2,2}[l]$ is a linear combination of brackets of the form
$[a,x_1]$ and $[b,x_2]$, the proposition is proved.
\end{proof}

\begin{theorem}
\label{2VarGen-thm}
For $l\ge2$, the $l-1$ elements of the form $[x_1^i,x_2^{l-i}]$ for $i = 1,\ldots,l-1$ constitute
a basis for $B_{2,2}[l]$, so for any $i,j\geq 1$, $B_{2,2}[(i,j)] = \mathbb{C}\cdot[x_1^i,x_2^j]$.
\end{theorem}

\begin{proof}
We will show that $\dim B_{2,2}[l]\geq l-1$.  Since we have already found $l-1$ generators for
$B_{2,2}[l]$, we conclude that $\dim B_{2,2}[l]$ must be equal to
$l-1$, and thus the spanning set we found must be
a basis for $B_{2,2}[l]$.

We claim that $[x_1^{l-1},x_2]$ is non-zero, i.e. $[x_1^{l-1},x_2]$ is not in $[[A_2,A_2],A_2]$. Note that $[[A_2,A_2],A_2]$ is spanned by elements of the form $[[m_1,m_2],m_3]$,
where $m_1$, $m_2$, and $m_3$ are monomials in $A_2$, and the only brackets of
this form which contain either $x_1^{l-1}x_2$ or $x_2 x_1^{l-1}$ are either of the form
$[[x_1^i,x_2],x_1^{l-1-i}]$ or $[[x_2,x_1^i],x_1^{l-1-i}] = -[[x_1^i,x_2],x_1^{l-1-i}]$.  In these
brackets, the coefficients of $x_1^{l-1}x_2$ and $x_2x_1^{l-1}$ are always equal.
Therefore, no linear combination of these brackets can give opposite signs on
$x_1^{l-1}x_2$ and $x_2x_1^{l-1}$, as in $[x_1^{l-1},x_2]$.  Hence, $[x_1^{l-1},x_2]$ is not in
$[[A_2,A_2],A_2]$.

Consider the Lie algebra $\mathfrak{gl}(2, \mathbb{C})$. Then it
has an action on $B_{2,2}[l]$ since it has a natural action on
the generators $\{x_1,x_2\}$. By direct computation, we can see
that $[x_1^{l-1},x_2]$ is a highest weight vector 
for $\mathfrak{gl}(2, \mathbb{C})$ with weight $(l-1,1)$. From
the representation theory of $\mathfrak{gl}(2,\mathbb{C})$, it
follows that this vector generates an $(l-1)$-dimensional
irreducible representation 
of $\mathfrak{gl}(2,\mathbb{C})$ contained in $B_{2,2}[l]$. 

Hence, $\dim B_{2,2}[l] \geq l-1$,
and the conclusion follows from the argument given in the beginning of the proof.
\end{proof}

\subsection{The $n=3$ Case}
\begin{proposition}
\label{3VarGen}
For $l\geq 2$, the $l^{2}-1$ non-zero elements of the form $[x_1^{i_{1}}
x_2^{i_{2}},x_3^{i_{3}}]$ and
$[x_1^{i_{1}} x_3^{i_{3}},x_2^{i_{2}}]$ for $i_{1}+i_{2}+i_{3}=l$ constitute a spanning set for $B_{3,2}[l]$.
\end{proposition}
\begin{proof}
By a similar argument as before, every element of $B_{3,2}$ can be expressed as a linear combination
of the brackets $[a,x_1]$, $[b,x_2]$, and $[c,x_3]$, where $a$, $b$, and $c$ are monomials with
degree no less than 1.

Consider $[a,x_1]$.  This may be written as a constant multiple of a bracket of the form $[x_1^{i_{1}},a_1]$,
where $a_1$ is a product of only $x_2$'s and $x_3$'s.  We then write $[x_1^{i_{1}},a_1]$ as the sum of
$[a_1 x_1^{i_{1}-1},x_1]$ with brackets of the form $[x_1^{i_{1}} d_1,x_2]$ and $[x_1^{i_{1}} d_2,x_3]$, where $d_1$
and $d_2$ are products of only $x_2$'s and $x_3$'s.

We can write $[x_1^{i_{1}} d_1,x_2] = k_1[x_1^{i_{1}} x_3^{i_{3}},x_2^{i_{2}}]$ and $[x_1^{i_{1}} d_2,x_3] = k_2[x_1^{i_{1}} x_2^{i_{2}},x_3^{i_{3}}]$.
Noting that $[a_1 x_1^{i-1},x_1]$ is a constant multiple of $[x_1^{i_{1}},a_1]$, we may now solve for $[x_1^{i_{1}},a_1]$,
realizing it as a linear combination of $[x_1^{i_{1}} x_3^{i_{3}},x_2^{i_{2}}]$ and $[x_1^{i_{1}} x_2^{i_{2}},x_3^{i_{3}}]$.  Therefore
$[a,x_1]$ is a linear combination of $[x_1^{i_{1}} x_3^{i_{3}},x_2^{i_{2}}]$ and $[x_1^{i_{1}} x_2^{i_{2}},x_3^{i_{3}}]$.

By performing a similar analysis on $[b,x_2]$ and $[c,x_3]$, we find that every element of $B_{3,2}$ can be
expressed as a linear combination of $[x_1^{i_{1}} x_3^{i_{3}},x_2^{i_{2}}]$, $[x_1^{i_{1}} x_2^{i_{2}},x_3^{i_{3}}]$, and $[x_2^{i_{2}} x_3^{i_{3}},x_1^{i_{1}}]$.
Noting that $[x_1^{i_{1}} x_3^{i_{3}},x_2^{i_{2}}] + [x_1^{i_{1}} x_2^{i_{2}},x_3^{i_{3}}] + [x_2^{i_{2}} x_3^{i_{3}},x_1^{i_{1}}] = 0$, we have that every
element of $B_{3,2}$ can be expressed as a linear combination of $[x_1^{i_{1}} x_3^{i_{3}},x_2^{i_{2}}]$ and
$[x_1^{i_{1}} x_2^{i_{2}},x_3^{i_{3}}]$.
\end{proof}

By using a similar method to the two variable case, 
we can consider the $\mathfrak{gl}(3,\mathbb{C})$ action 
and find that the element $[x_1^{l-1},x_2]$ (which we showed to
be nonzero when considering the two variable case) is a highest
weight vector of weight $(n-1,1,0)$. It then follows from the
representation theory of $\mathfrak{gl}(3,\mathbb{C})$ that this
vector generates a representation of dimension $l^2-1$, and hence
the dimension of $B_{3,2}[l]$ is at least $l^{2}-1$.
Combining with Proposition~\ref{3VarGen} we have

\begin{theorem}\label{3VarGen-thm}
For any ${\mathbf i}=(i_1,i_2,i_3)\in (\mathbb{Z}_{>0})^{3}$, $[x_1^{i_1} x_2^{i_2},x_3^{i_3}]$ and $[x_1^{i_1} x_3^{i_3},x_2^{i_2}]$
constitute a basis for $B_{3,2}[{\mathbf i}]$.
\end{theorem}

\begin{remark}
The proof of Proposition~\ref{2VarGen-thm}  mimics the
manipulation of 1-forms done by 
Feigin and Shoikhet~\cite{FS},
except in the language of brackets. 

Theorems~\ref{2VarGen-thm}
and \ref{3VarGen-thm} can
be obtained from Theorems 1.3 and 1.4 in \cite{FS}, respectively. But here we 
have given new direct proofs of these theorems without using the
results of \cite{FS}. We have not been able to generalize these
proofs to the case $n\ge 4$.  
\end{remark}

\section{The structure of $B_{n,2}$ for general $n$}

\subsection{The Main Theorem about $B_{n,2}$}
Let $P_n$ be the set of all permutations of ${1,2,\ldots ,n}$ which have the form
$(2, 3)^{\delta_{2}} (3, 4)^{\delta_{3}} \cdots 
(n-1, n)^{\delta_{n-1}}$, where $(i, j)$ is the permutation of $i$ and $j$, and $\delta_i = 0$ or $1$. From now on, let ${\mathbf i}=(i_1,\ldots ,i_n)\in(\mathbb{Z}_{>0})^{n}$.

\begin{theorem}
\label{Basis}
The basis elements for $B_{n,2}[\bf{i}]$ are the $2^{n-2}$
brackets given by \linebreak $[x_{p(1)}^{i_{p(1)}}\cdots x_{p(n-1)}^{i_{p(n-1)}},x_{p(n)}^{i_{p(n)}}]$ for
$p\in P_n$. 

In particular, we have $$\dim B_{n,2}[{\mathbf i}] = 2^{n-2}.$$
\end{theorem}

\begin{remark}\label{nonre}
Note that if some $i_s$ is zero then a basis of
$B_{n,2}[{\mathbf i}]$ is given by Theorem \ref{Basis} 
for a smaller number of variables. Thus, Theorem \ref{Basis} provides a basis of $B_{n,2}[l]$ for
any $l$, and thus a homogeneous basis of $B_{n,2}$. An
interesting property of this basis is that it consists of
elements whose monomials are non-redundant 
(i.e. every letter occurs only once in some power).
\end{remark}

The proof of Theorem \ref{Basis} is given in the next three
subsections. 

\subsection{The Feigin-Shoikhet Isomorphism}

We will use the isomorphism in \cite{FS} between
$B_{n,2}$ and $\Omega^{\mathrm{even}+}_{\mathrm{closed}}(\mathbb{C}^n)$, the closed even differential forms with positive
degree, to prove the main theorem. 
Recall $\phi_n$ is a homomorphism of algebras: 
$$
\phi_n: A_n \rightarrow \Omega^{\mathrm{even}}(\mathbb{C}^n)_{\ast}
$$ 
which takes $x_i \in A_n$
to $x_i \in \Omega^0(\mathbb{C}^n)$ and
$$
\phi_n(x_i x_j) = x_i \ast x_j = x_i x_j + dx_i\wedge dx_j.
$$
Feigin and Shoikhet proved that $\phi_n$ induces an isomorphism 
$$
\phi_n: B_{n,2} \rightarrow
\Omega^{\mathrm{even}+}_{\mathrm{closed}}(\mathbb{C}^n).
$$

Now let $w \in
\Omega^p(\mathbb{C}^n)$ be a $p$-form.  We say that $w$ has multidegree ${\mathbf i}$ if
every $x_s$ occurs $i_s$ times in every monomial of $w$.  Define $\Omega^p(\mathbb{C}^n)[{\mathbf i}]$
to be the space of all forms of multidegree ${\mathbf i}$.

The main theorem is then a consequence of the following two lemmas:

\begin{lemma}
\label{NumGen}
$\dim \Omega^{\mathrm{even}+}_{\mathrm{closed}}(\mathbb{C}^n)[{\mathbf i}] = 2^{n-2}$.
\end{lemma}

\begin{lemma}
\label{IndGen}
The $2^{n-2}$ brackets described in the main theorem are
linearly independent.
\end{lemma}

\subsection{Proof of Lemma~\ref{NumGen}}

We first prove a more basic lemma, from which Lemma~\ref{NumGen} will follow.

\begin{lemma}
$\dim \Omega^p_{\mathrm{closed}}(\mathbb{C}^n)[{\mathbf i}] = \dbinom{n-1}{p-1}$.
\end{lemma}

\begin{proof}
By the Poincar\'e Lemma, the De Rham differential defines an
isomorphism $d:
\Omega^{p-1}(\mathbb{C}^n)[{\mathbf i}]/\Omega^{p-1}_{\mathrm{closed}}(\mathbb{C}^n)[{\mathbf i}]\to
\Omega^p_{\mathrm{closed}}(\mathbb{C}^n)[{\mathbf i}]$.

Hence, if $D(p):=\dim \Omega^{p}_{\mathrm{closed}}(\mathbb{C}^n)[{\mathbf i}]$, we have
the recurrence relation:
\begin{equation*}
D(p) = \dbinom{n}{p-1} - D(p-1),\text{ and }D(0) = 0.
\end{equation*}
A simple inductive argument shows $D(p) = \dbinom{n-1}{p-1}$, as desired.
\end{proof}
Lemma~\ref{NumGen} now follows from a simple combinatorial identity:
$$
\dim \Omega^{\mathrm{even}+}_{\mathrm{closed}}(\mathbb{C}^n)[{\mathbf i}] 
=\displaystyle \sum_{k=1}^{\infty}{\dim \Omega^{2k}_{\mathrm{closed}}(\mathbb{C}^n)[{\mathbf i}]} 
= \displaystyle \sum_{k=1}^{\infty}{\dbinom{n-1}{2k-1}} = 2^{n-2}.
$$

\subsection{Proof of Lemma~\ref{IndGen}}
We begin by computing the image under the map $\phi_n$ of the brackets with
the form given in the statement of the main theorem.

\begin{lemma}
\label{PhiProd}
We have
$$
\phi_n(x_1^{i_1}\ldots x_n^{i_n}) = x_1^{i_1}\ldots x_n^{i_n} \displaystyle \sum_{\substack{S\subset
\{1,\ldots ,n\}\\|S| \textrm{ even}}}{\bigwedge_{k\in S} i_k \frac{dx_k}{x_k}},
$$
where the indices in the wedge product are in increasing order.
\end{lemma}

\begin{proof}
We prove this by induction.  For $n=1$, we have $\phi_n(x_1^{i_1}) = x_1^{i_1}$, as desired.
Assume the lemma is true for $n$.  Then 
$$
\phi_n(x_1^{i_1}\cdots x_{n+1}^{i_{n+1}}) =
\phi_n(x_1^{i_1}\cdots x_n^{i_n})\ast x_{n+1}^{i_{n+1}} = \left(x_1^{i_1}\cdots x_n^{i_n}
\displaystyle\sum_{\substack{S\subset\{1,\ldots ,n\}\\|S| \textrm{ even}}}{\bigwedge_{k\in S}
i_k \frac{dx_k}{x_k}}\right)\ast x_{n+1}^{i_{n+1}}.
$$

Note that in the expansion of the last expression, the sum of the $2l$-forms comes from
$d[(2l-2)$-forms$]\wedge dx_{n+1}^{i_{n+1}} + 2l$-forms $\wedge x_{n+1}^{i_{n+1}}$.
The first term gives all $2l$-forms which contain $dx_{n+1}$, whereas the second term gives
all $2l$-forms which do not contain $dx_{n+1}$.  Together, all possible $2l$-forms appear
in the expansion.  These forms correspond to the subsets $S$ of $\{1,\ldots ,n+1\}$ with exactly
$2l$ elements.  It is not hard to see that the coefficients of these forms are precisely
the ones in the lemma.
\end{proof}

By direct computation, we have:

\begin{corollary}
\label{PhiBrac}
Let $\omega_{S}(x_{1},\ldots,x_{n})=2x_1^{i_1}\cdots x_n^{i_n}\bigwedge_{k\in S} i_k
\frac{dx_k}{x_k}$ for $S\subset \{1,\ldots, n\}$.
Then for $i_1,\ldots ,i_n > 0$, we have
$$
\phi_n([x_1^{i_1}\cdots x_{n-1}^{i_{n-1}},x_n^{i_n}]) =
\sum_{\substack{n\in S\subset \{1,\ldots ,n\}\\|S| \textrm{ even}}}\omega_{S}(x_{1},\ldots,x_{n}),
$$
where the indices in the wedge product are in increasing order.
\end{corollary}

Notice that for $p\in P_{n}$, we have 
$$\phi_n([x_{p(1)}^{i_{p(1)}}\cdots x_{p(n-1)}^{i_{p(n-1)}},x_{p(n)}^{i_{p(n)}}]) =
\sum_{\substack{p(n)\in S\subset \{1,\ldots ,n\}\\|S| \textrm{ even}}}\epsilon(S,p)\ \omega_{S}(x_{p(1)},\ldots,x_{p(n)}),$$ where $\epsilon(S,p)=\pm1$, depending on the choice of $S$ and $p$.

Denote $\omega_{S}(x_{p(1)},\ldots,x_{p(n)})$ by
$\omega_{S}^{p}$. We are now ready to prove 
Lemma~\ref{IndGen}.

\begin{proof}[Proof of lemma~\ref{IndGen}]
We proceed by induction on the number of variables.  It is easy to see that the lemma is true for $n=2,3$ from the results in Section~\ref{n23}.
Assume the lemma is true up to $n\geq 3$.  We now prove the lemma for $n+1$ variables.

Let $P_{n+1}^1$ be the set
of permutations which contain $(n, n+1)$ (i.e., with
$\delta_n=1$), and $P_{n+1}^{2}$ be its complement in $P_{n+1}$. Then we have 
$P_{n+1} = P_{n+1}^1 \cup P_{n+1}^2$.

By applying the isomorphism $\phi_{n+1}$, it is enough to show that the $2^{n-1}$ forms:
\begin{equation*}
\omega^{p}=\sum_{\substack{p(n+1)\in S\subset \{1,\ldots ,n+1\}\\|S| \textrm{ even}}}\epsilon(S,p)\ \omega_{S}^{p}
\end{equation*}
are linearly independent. 

For any $p\in P^{1}_{n+1}$, the components of $\omega^{p}$ which do not contain $dx_{n+1}$ are precisely the forms in the $n$ variable case
which appear in $\omega^{p'}$, where $p'\circ (n,n+1) = p$.
Hence, the $\omega^{p}$ for $p\in P^{1}_{n+1}$ are linearly independent.

Furthermore, since every form which appears in $\omega_{S}^{p}$ for $p\in P^{2}_{n+1}$ contains $dx_{n+1}$,
we only need to show that the forms $\omega^{p}$ with $p\in P_{n+1}^{2}$ are linearly independent.

Let $\mathcal{S}=\{S\subset \{1,\ldots,n+1\}|1,n+1\in S \text{ and } |S| \text{ is even}\}$. 
For any $p\in P^{2}_{n+1}$, $\displaystyle\sum_{S\in \mathcal{S}}\epsilon(S,p)\omega_{S}^{p}$
is a linear combination of even forms containing $dx_{1}\wedge dx_{n+1}$. It is enough to show that these $2^{n-2}$ sums
are linearly independent.

It suffices to prove the invertibility of the $2^{n-2}\times 2^{n-2}$ matrix where each row represents
a bracket $p \in P^{2}_{n+1}$, each column represents a form $S \in \mathcal{S}$, and whose entries are the $\epsilon(S,p)$'s.
For the rows, we choose
the order recursively, beginning with the identity permutation.
Given the first $2^k$ elements, the next $2^k$ elements are given by composition with $(k+2, k+3)$.

For the columns, we will represent the form $dx_{j_1}\wedge\cdots \wedge dx_{j_m}$ by the ordered
$m$-tuple $(j_1,\ldots ,j_m)$.  We again choose the
order recursively, beginning with $(1,n+1)$. Given the first $2^k$ columns, the next $2^k$ columns are given
by appending $k+2,k+3$ to the first $2^{k-1}$ columns
and by replacing $k+2$ with $k+3$ in the next $2^{k-1}$ columns.
 
We prove the invertibility of this matrix by induction on $n$.
When $n=3$, the matrix is given by
$\left( \begin{array}{cc}
1 & 1 \\
1 & -1 \end{array} \right)$,
which is clearly invertible.  Assume it is true for $n\geq 3$.

Divide the matrix into equal fourths.  Call the submatrices $\alpha_n$, $\beta_n$, $\gamma_n$, $\delta_n$.
Note that $\alpha_n$ is the matrix for the $n$ variable case.  Now further divide each of these
submatrices into four more equal quadrants.  Call them $\alpha_n^{1,1}$, $\alpha_n^{1,2}$, $\alpha_n^{2,1}$,
$\alpha_n^{2,2}$, etc.  In the case of $\alpha_n$, we have $\alpha_n^{1,1}$ = $\alpha_{n-1}$,
$\alpha_n^{1,2}$ = $\beta_{n-1}$, $\alpha_n^{2,1}$ = $\gamma_{n-1}$, $\alpha_n^{2,2}$ = $\delta_{n-1}$.

Because changing the
position of $n$ in the permutation has no effect on the sign of the forms which do not contain
$n$, we have $\alpha_n = \gamma_n$ (and $\alpha_{n-1} = \gamma_{n-1}$). We also have $\alpha_n^{1,1} = \beta_n^{1,1}$ and $\alpha_n^{1,2} = \beta_n^{1,2}$ because the permutations in those
rows leave $n-1$ and $n$ fixed.  By similar analysis of the permutations, we can show the matrix has the form:
\begin{equation*}
\left(\begin{array}{cc}\alpha_n & \beta_n \\\gamma_n & \delta_n\end{array}\right)=\left( \begin{array}{cccc}
\alpha_{n-1} & \beta_{n-1} & \alpha_{n-1} & \beta_{n-1} \\
\alpha_{n-1} & \delta_{n-1} & * & \beta_{n-1} \\
\alpha_{n-1} & \beta_{n-1} & -\alpha_{n-1} & \beta_{n-1} \\
\alpha_{n-1} & \delta_{n-1} & * & \delta_{n-1} \end{array} \right).
\end{equation*}

Subtracting the last $2^{n-3}$ rows from the first $2^{n-3}$ rows gives
$$\left( \begin{array}{cccc}
\alpha_{n-1} & \beta_{n-1} & \alpha_{n-1} & \beta_{n-1} \\
\alpha_{n-1} &\delta_{n-1} & * & \beta_{n-1} \\
0 & 0 & 2\alpha_{n-1} & 0 \\
0 & 0 & * & \beta_{n-1}-\delta_{n-1} \end{array} \right).$$

It remains to show that $\alpha_n$, $\alpha_{n-1}$, and $\beta_{n-1}-\delta_{n-1}$ are invertible.
$\alpha_n$ and $\alpha_{n-1}$ are invertible by the induction hypothesis.  We see that
$\beta_{n-1}-\delta_{n-1}$ is invertible by subtracting the last half of the rows from
the first half in the invertible matrix
$$
\alpha_n =
\left( \begin{array}{cc}
\alpha_{n-1} & \beta_{n-1} \\
\alpha_{n-1} & \delta_{n-1} \end{array} \right).
$$
\end{proof}

\subsection{Finite order case}

\begin{theorem}\label{basis-rel}
For any ${\mathbf i}= (i_1,\ldots ,i_n)\in (\mathbb{Z}_{>0})^{n}$,
\begin{equation*}
B_{n,2}^R[{\mathbf i}]= \left\{ \begin{array}{ll}
        0 & \mbox{if $i_s \geq m_s$ for some $s$};\\
        B_{n,2}[{\mathbf i}] & \mbox{if $i_s < m_s$ for all $s$}.\end{array} \right.
\end{equation*}
\end{theorem}

\begin{proof}
It is clear that if $i_s<m_s$ for all $s$ then the relations have
no effect, so the statement of the theorem holds. 
Now assume that for some $s$, $i_s\ge m_s$. 
Then the images in $B_{n,2}^R$ of all the basis elements from 
Theorem \ref{Basis} are zero. But these elements must span
$B_{n,2}^R[{\mathbf i}]$, which implies that this space is
zero, as desired. 
\end{proof}

\begin{remark}
In this proof it is important that the basis elements involve
only non-redundant monomials, see remark \ref{nonre}.
\end{remark}

\section{The structures of $B_{2,m}[r,1]$ and $B_{2,m}[r,2]$}

Let $W_{n}$ be the Lie algebra of polynomial vector fields on $\mathbb{C}^{n}$. In \cite{FS}, Feigin and Shoikhet described an action of $W_n$ on $B_{n,k}$. 

From now on, let $A=A_{2}$ be the free algebra generated by $x,y$
and $L_{i}=L_{i}(A)$.
We denote by  $L_{i}[r,s]$ (and $B_{2,m}[r,s]$) the space of elements of $L_{i}$ (and $B_{2,m}$) with multi-degree $(r,s)$, that is, consisting of monomials having $r$ copies of $x$ and $s$ copies of $y$. The purpose of this section is to compute the bases of $B_{2,m}[r,1]$ and $B_{2,m}[r,2]$ which will help us to find the $W_2$-module structures of $B_{2,3}$ and $B_{2,4}$, and obtain some information about the structure of $B_{2,m}$ for general $m$ in the subsequent sections. 

Define $\mathrm{ad}_{a}b=[a,b]$. Then we introduce the following elements:
\begin{eqnarray*}
b^{(l)}_{i,j,k}&=&\mathrm{ad}^{i}_{x}\circ \mathrm{ad}_{y}\circ \mathrm{ad}^{j}_{x}\circ \mathrm{ad}_{x^{k}}(y^{l});\\
b^{(l)}_{i,j}&=&\mathrm{ad}^{i}_{x}\circ \mathrm{ad}_{x^{j}}(y^{l}).
\end{eqnarray*}

Notice that $b^{(l)}_{i,j,k}$ is an element in $L_{i+j+3}[i+j+k,l+1]$, and $b^{(l)}_{i,j}$ is an element in $L_{i+2}[i+j,l]$. For simplicity, when $l=1$, we denote $b^{(1)}_{i,j,k}$ by $b_{i,j,k}$ and $b^{(1)}_{i,j}$ by $b_{i,j}$.

\subsection{Structure of $B_{2,m}[r,1]$}

\begin{theorem}\label{B_{2,m}[r,1]}  
For $m \geq 2$ we have:
$$B_{2,m}[r,1]  = \left\{ \begin{array}{ll}
0, & \textrm{$r\leq m-2;$}\\
\mathbb{C}\cdot b_{m-2,r-m+2},& \textrm{$r\geq m-1.$}
\end{array} \right.$$
\end{theorem}

First we prove two lemmas.

\begin{lemma}\label{d/dx surjective} 
For $r \geq 1$ and $s \geq 0$, the linear map $\frac{\partial}{\partial x}:\ A_2[r,s] \rightarrow A_2[r-1,s]$ is surjective. 
\end{lemma}

\begin{proof} 
We do induction on $s$. For $s=0$ the statement is obviously true. 
Now suppose $s>0$. Every monomial in $A_2[r-1,s]$ has the form $ m y x^a $, where $0 \leq a \leq r-1$ and $m$ is a monomial in $A_2[r-a-1,s-1]$. By the induction hypothesis there exists a polynomial $p$ such that $\frac{\partial}{\partial x}p = m$. 

Now we show  by induction on $a$ that there exists a polynomial $q$ such that $\frac{\partial}{\partial x}q = m y x^a$. If $a=0$, put $q=py$. Suppose $a>0$, then by the induction hypothesis there exists a polynomial $f$ in $A_2[r,s]$ such that $\frac{\partial}{\partial x}f = p y x^{a-1}$. Then $q= pyx^{a} - af$ will be a solution. 
\end{proof}

\begin{lemma}\label{ker d/dx [r,1]} The kernel of the map $\frac{\partial}{\partial x}:\ A[r,1] \rightarrow A[r-1,1]$ is $\mathbb{C}\cdot b_{r-1,1}$.
\end{lemma}

\begin{proof} 
By Lemma \ref{d/dx surjective}, $\dim\ker \frac{\partial}{\partial x} = \dim A[r,1] - \dim A[r-1,1] = (r+1)-r=1$. The element $b_{r-1,1}$ is in the kernel of $\frac{\partial}{\partial x}$, and it is non-zero. Thus the lemma is proved.
\end{proof}

Now we can prove the  theorem.

\begin{proof}[Proof of Theorem \ref{B_{2,m}[r,1]}] 
We first prove that for $m \leq r+1$ the element $b_{m-2,r-m+2}$ spans $B_{2,m}[r,1]$. 
We do this by induction on $r$. For $r=1$, $[x,y]$ spans $B_{2,2}[1,1]$. 
Now suppose $r>1$. 
The statement is true for $m=r+1$ obviously. Suppose $m \leq r$. 
For any $w\in B_{2,m}[r,1]$, we have $\frac{\partial}{\partial x} w\in B_{2,m}[r-1,1]$. By induction hypothesis\linebreak $\frac{\partial}{\partial x} w = cb_{m-2,r-m+1}$ for some constant $c$. 
Let $p=\frac{c}{r-m+2} b_{m-2,r-m+2}$. By Lemma \ref{ker d/dx [r,1]} we have $w-p \in \ker \frac{\partial}{\partial x} \subseteq L_{r+1}$. 
Since $m \leq r$, $w=p\in B_{2,m}$. So the statement is proved by induction.

Therefore we have $\dim B_{2,m}[r,1] \leq 1$ for $m \leq r+1$, and  $\dim B_{2,m}[r,1]=0$ for $m>r+1$. Since $\sum_{1 \leq m \leq r+1}\dim B_{2,m}[r,1]=\dim A_2[r,1]=r+1$ and $\dim B_{2,1}[r,1]=1$, we have $\dim B_{2,m}[r,1] = 1$ for $m \leq r+1$.     
\end{proof}

\subsection{Structure of $B_{2,m}[r,2]$}

\begin{theorem} \label{basis of A_2[r,2]} 
For $m \geq 2$ we have:
$$\dim B_{2,m}[r,2]  = \left\{ \begin{array}{ll}
m-1, & \textrm{$m\leq r+1;$}\\
\lfloor \frac{r+1}{2} \rfloor, & \textrm{$m=r+2.$}
\end{array} \right.$$
A basis of $B_{2,m}[r,2]$ for $ m \leq r+1$ is given by the $m-1$ elements 
$$b_{i,j,r-m+3} \text{ for } i+j=m-3, \text{ and } b^{(2)}_{m-2,r-m+2}.$$ 
\end{theorem}         

Before starting the proof we will prove several lemmas.



\begin{lemma} \label{basis of B_{2,r+2}} 
The set 
$S_{r}=\{b_{i,j,1}|i+j=r-1, j \text{ is even } \}$ is a basis of $L_{r+2}[r,2]$.
\end{lemma}

\begin{proof} 
At first, we prove elements in $S_{r}$ are independent by induction on $r$. For \linebreak $r=1$ the claim is obvious. Assume it is true for $r-1$. 
If $r$ is even, then these elements have the form $[x,b_{i,j,1}]$ where $i+j=r-2$ and $j$ is even. These elements are independent by the induction hypothesis because $\sum {\alpha}_{i,j} [x,b_{i,j,1}] =[x, \sum {\alpha}_{i,j} b_{i,j,1}]$ has the leading monomial $xm$ where $m$ is the leading monomial of $\sum {\alpha}_{i,j} b_{i,j,1}$.

If $r$ is odd, by a similar argument as the even case, we only need to show that the element $b_{0,r-1,1}$ is independent from the others. Since it is the only element which has the monomial $2yx^ry$, the conjectured basis elements are independent by induction. 

Now we show that elements in $S_{r}$ span $L_{r+2}[r,2]$. It is enough to show that $b_{0,r-1,1}$ with even $r$ is in $[x,L_{r+1}]$. Applying Jacobi identity repeatedly, we obtain
\begin{eqnarray*}
b_{0,r-1,1}&=& [[[y,x],x],[x,\ldots[x,y]\ldots]]+[x,L_{r+1}]\\
&=&\ldots= [[[y,x],\ldots],x],[x,\ldots[x,y]\ldots]]+[x,L_{r+1}],
\end{eqnarray*}
where the last element has equal number $i$ of copies of $x$ in the first and the second major brackets. But this element is zero, so the original element $b_{0,r-1,1}$ is in $[x,L_{r+1}]$. 
\end{proof}

\begin{lemma} \label{basis of B_{2,r+1}} 
The set $S_{r}'=S_{r}\cup \{b_{i,j,2}|i+j=r-2\}\cup\{b^{(2)}_{r-1,1}\}$ is a basis of $L_{r+1}[r,2]$. In particular, the set $S_{r}'-S_{r}$ is a basis of $B_{2,r+1}[r,2]$.
\end{lemma}

\begin{proof} At first, we prove by induction on $r$ that the elements in $S'_{r}$ are independent. When $r=1$, it is easy to see. Now suppose $r>1$. 

If $r$ is even, all the elements in $S'_{r}$ except the element $b_{0,r-2,2}$ will have the form $[x,b]$ where $b$ is in  $S_{r-1}'$ which is the basis of $L_{r}[r-1,2]$. As in the proof of Lemma \ref{basis of B_{2,r+2}}, we only need to show $b_{0,r-2,2}$ is independent from the others. Since its  leading monomial is $2yx^ry$, which is not found in the others, the elements in $S'_{r}$ are independent.  

If $r$ is odd, the elements in $S'_{r}$ are $[x,b]$ where $b\in S'_{r-1}$, and two other elements $b_{0,r-2,2}$, $b_{0,r-1,1}$. We observe that $b_{0,r-1,1}$ has a leading monomial $2yx^ry$ which no other elements in $S'_{r}$ have, therefore $b_{0,r-1,1}$ is independent from them. 

Now let $r=2i+1$. By direct computation, we have 
\begin{eqnarray*}
b_{0,r-2,2}&=&(2i-1)xyx^{2i}+(-2i^2+3i)x^2yx^{2i-1}y+0xyx^{2i-1}yx+\cdots,\\
b_{1,r-3,2}&=&2xyx^{2i}+(-2i+2)x^2yx^{2i-1}y+0xyx^{2i-1}yx+\cdots,\\ 
b_{2,r-3,1}&=&0xyx^{2i}+2x^2yx^{2i-1}y-4xyx^{2i-1}yx+\cdots.
\end{eqnarray*}
Since these monomials are not present in the other elements of $S'_{r}$, we have that the element $b_{0,r-2,2}$ is independent from the other elements. Therefore by induction all the elements of $S'_{r}$ are independent.  

Now we show that $S'_{r}$ is a spanning set by induction on $r$. For $r=1$ the statement is true. Assuming the statement for $r-1$, we obtain that 
$[x,L_r[r-1,2]]$ is spanned by elements $[x,b]$ for $b\in S'_{r-1}$. Since the space $[y,L_r[r,1]]$ is spanned by $b_{0,r-1,1}$ we only need to show that  $[x^2,b_{0,r-3,1}]$ is in the spanning space of $S'_{r}$. By repeatedly applying Jacobi identity, we have:
\begin{eqnarray*}
[x^2,b_{0,r-3,1}]&=& [[x^2,y],[x,\ldots[x,y]\ldots]]+[y,L_r[r,1]]\\
& =& [[[x^2,y],x],[x,\ldots[x,y]\ldots]]+[x,L_r[r-1,2]]+[y,L_r[r,1]]\\
&=& \ldots = [\ldots[x^2,y],x],\ldots],x],y]+[x,L_r[r-1,2]]+[y,L_r[r,1]].
\end{eqnarray*}
So the set $S'_{r}$ spans $L_{r+1}[r,2]$ and we proved the lemma.
\end{proof}

\begin{lemma}\label{ker d/dx} 
The linear map $\frac{\partial}{\partial x}: A_2[r,2] \rightarrow A_2[r-1,2]$ has the property 
$$\ker \frac{\partial}{\partial x} \cap [A,A] \subseteq L_{r+1}[r,2].$$ 
\end{lemma}

\begin{proof}
From Lemma \ref{basis of B_{2,r+2}} 
and \ref{basis of B_{2,r+1}} it follows that 
$$\dim B_{2,r+1}[r,2]=r, \dim B_{2,r+2}[r,2]=\lfloor \frac{r+1}{2} \rfloor.$$
We have the induced linear map $\frac{\partial}{\partial x}|_{B_{2,r+1}[r,2]}: B_{2,r+1}[r,2] \rightarrow B_{2,r+1}[r-1,2]$ which is surjective because $\frac{\partial}{\partial x}b_{i,j,2}=2b_{i,j,1}$. So 
$$\dim \ker \frac{\partial}{\partial x}|_{B_{2,r+1}[r,2]}= \dim B_{2,r+1}[r,2] - \dim B_{2,r+1}[r-1,2] = r -\lfloor \frac{r}{2} \rfloor = \lfloor \frac{r+1}{2} \rfloor.$$ 
Also $\frac{\partial}{\partial x}$ maps $B_{2,r+2}[r,2]$ to zero, so $\dim(\ker \frac{\partial}{\partial x} \cap L_{r+1}[r,2])=2 \lfloor \frac{r+1}{2} \rfloor$.
 
If $r$ is odd, $2 \lfloor \frac{r+1}{2} \rfloor = r+1$, so $\ker \frac{\partial}{\partial x} \subseteq L_{r+1}[r,2]$. If $r$ is even, $2 \lfloor \frac{r+1}{2} \rfloor = r$. In this case we consider the induced map $\frac{\partial}{\partial x}|_{B_{2,1}[r,2]}: B_{2,1}[r,2] \rightarrow B_{2,1}[r-1,2]$. These spaces are the spaces of cyclic words, so $\dim \ker \frac{\partial}{\partial x}|_{B_{2,1}[r,2]} \geq \dim B_{2,1}[r,2] - \dim B_{2,1}[r-1,2] = \lceil \frac{r+1}{2} \rceil - \lceil \frac{r}{2} \rceil = 1$ if $r$ is even. So for even $r$, a one-dimensional subspace of $\ker \frac{\partial}{\partial x}$ lies in $B_{2,1}$. Therefore $\ker \frac{\partial}{\partial x} \cap [A,A] \subseteq L_{r+1}[r,2]$.
\end{proof}

Now we prove Theorem \ref{basis of A_2[r,2]}.

\begin{proof}[Proof of Theorem \ref{basis of A_2[r,2]}]

Let us prove that any element $w$ of $B_{2,m}[r,2]$ is a linear combination of the conjectured basis elements. We do induction on $r$. 

If $r=1$, both $[x,y^2]$ and $[y,[x,y]]$ are basis elements. 
If $r>1$, we have $\frac{\partial}{\partial x} w$ which has degree $r-1$ in $x$. So by induction hypothesis $\frac{\partial}{\partial x} w = \sum_{i+j=m-3} \alpha_{i,j} b_{i,j,r-m+2}+ \alpha b^{(2)}_{m-2,r-m+1}$. Put 
$p=\sum_{i+j=m-3} \frac{\alpha_{i,j}}{r-m+3}b_{i,j,r-m+3}+ \frac{\alpha}{r-m+2}b^{(2)}_{m-2,r-m+2}$, then $\frac{\partial}{\partial x}(w-p)=0$.
So $w-p \in \ker \frac{\partial}{\partial x}\cap [A,A] \subseteq L_{r+1}$ by Lemma \ref{ker d/dx}. So $w=p$ in $L_m/L_{m+1}$ $(m \leq r)$, and $p$ is a combination of basis elements. Therefore the required elements span $B_{2,m}[r,2]$ and $\dim B_{2,m} \leq m-1$ ($2 \leq m \leq r$). 

We know that $\dim B_{2,1}=\lceil \frac{r+1}{2} \rceil$, $\dim B_{2,r+2}=\lfloor \frac{r+1}{2} \rfloor$ and $\dim B_{2,m} \leq m-1$ ($2 \leq m \leq r$). But these numbers have to sum to $\dim A_2[r,2] = \frac{(r+1)(r+2)}{2}$, so $\dim B_{2,m} = m-1$ ($2 \leq m \leq r$) and the found spanning elements actually form a basis for $B_{2,m}[r,2]$.
 \end{proof}

\section{The multiplicities of $\mathcal{F}_{(p,1)}$ and $\mathcal{F}_{(p,2)}$ in $B_{2,m}$}

We consider the $W_{n}$-modules  on which the Euler vector field \linebreak $e=\sum_{i=1}^n x_i \frac {\partial}{\partial x_i}$ is semisimple with finite-dimensional eigenspaces and has its eigenvalues bounded from below. Let $W_n^0$ be the subalgebra of $W_n$ of vector fields vanishing at the origin.

Let ${\mathcal{F}}_D=\Hom_{U(W_n^0)}(U(W_n),F_D)$ be the irreducible $W_{n}$-module coinduced from a $\mathfrak{gl}(n,\mathbb{C})$-module $F_D$ where $D$ is a Young diagram having more than one column. (For reference about modules ${\mathcal{F}}_D$ see \cite{FF} or \cite{F}; for reference about Schur modules $F_D$ see \cite{Ful}). Let $(p,k)$, where $p \geq k$ are positive integers, denote a two-row Young diagram with $p$ boxes in the first row and $k$ boxes in the second row.

In this section we prove

\begin{theorem}\label{multiplicities} 
For $m \geq 3$, the $W_2$-module $B_{2,m}$ has in its Jordan-H\"older series one copy of the module $\mathcal{F}_{(m-1,1)}$, $\lfloor \frac{m-2}{2} \rfloor$ copies of $\mathcal{F}_{(m-1,2)}$, and $\lfloor \frac{m-3}{2} \rfloor$ copies of $\mathcal{F}_{(m-2,2)}$. The rest of the irreducible $W_2$-modules in the Jordan-H\"older series of $B_{2,m}$ are of the form $\mathcal{F}_{(p,k)}$ where $k \geq 3$. 
\end{theorem}

\begin{proof} 
If $B_{2,m}$ contains a module $\mathcal{F}_{(p)}$, then $\dim \mathcal{F}_{(p)}[r,0]= 1$ which contradicts $\dim B_{2,m}[r,0]=0$. Similarly, $B_{2,m}$ cannot contain the module of exact one-forms. Therefore all the irreducible $W_2$-modules contained in $B_{2,m}$ are of the form $\mathcal{F}_{(p,k)}$ where $k \geq 1$.   
  
At first, we find the multiplicities of the modules $\mathcal{F}_{(p,1)}$ in $B_{2,m}$. Notice that for modules $\mathcal{F}_{(p,k)}$ where $k\geq 2$, we have $\mathcal{F}_{(p,k)}[r,1]=0$. We also have 
$$ 
\dim \mathcal{F}_{(p,1)}[r,1]  = \left\{ \begin{array}{ll}
0, & \textrm{$r\leq p-1;$}\\
1, & \textrm{$r\geq p.$}
\end{array} \right.
$$
Comparing this to Theorem \ref{B_{2,m}[r,1]}, we obtain that $B_{2,m}$ has one copy of $\mathcal{F}_{(m-1,1)}$ and none of the other modules $\mathcal{F}_{(p,1)}$ where $p \neq m-1$. 

Now let us find the multiplicities of the modules $\mathcal{F}_{(p,2)}$ in $B_{2,m}$. For modules $\mathcal{F}_{(p,k)}$ where $k\geq 3$ we have $\mathcal{F}_{(p,k)}[r,2]=0$. We notice that
\begin{equation} \label{F_(p,2)[r,2]} 
\dim \mathcal{F}_{(p,2)}[r,2]  = \left\{ \begin{array}{ll}
0, & \textrm{$r\leq p-1;$}\\
1, & \textrm{$r\geq p.$}
\end{array} \right.
\end{equation}
We also have
$$
\dim \mathcal{F}_{(m-1,1)}[r,2]  = \left\{ \begin{array}{ll}
0, & \textrm{$r\leq m-3;$}\\
1, & \textrm{$r=    m-2;$}\\
2, & \textrm{$r\geq m-1.$}
\end{array} \right.
$$
By Theorem \ref{basis of A_2[r,2]}  we have
$$
\dim B_{2,m}[r,2] - \dim \mathcal{F}_{(m-1,1)}[r,2] = \left\{ \begin{array}{ll}
0, & \textrm{$r\leq m-3;$}\\
\lfloor \frac{m-3}{2} \rfloor, & \textrm{$r= m-2;$}\\
m-3, & \textrm{$r\geq m-1.$}
\end{array} \right.
$$
From formula (\ref{F_(p,2)[r,2]}), we have that $B_{2,m}$ has $\lfloor \frac{m-3}{2} \rfloor$ copies of $\mathcal{F}_{(m-2,2)}$ and $\lfloor \frac{m-2}{2} \rfloor$ copies of $\mathcal{F}_{(m-1,2)}$, which together with the module $\mathcal{F}_{(m-1,1)}$ account for the dimensions of $B_{2,m}[r,2]$. Finally, we remark that there may be some copies of the modules $\mathcal{F}_{(p,k)}$ with $k\geq 3$ in $B_{2,m}$ which we cannot detect with the help of the structures of $B_{2,m}[r,1]$ and $B_{2,m}[r,2]$.     
\end{proof}

We make the statement of this theorem more precise with the following

\begin{proposition} The module $\mathcal{F}_{(m-1,1)}$ is the last term of the Jordan-H\"older series of $B_{2,m}$, i.e. there is a projection map $B_{2,m} \twoheadrightarrow \mathcal{F}_{(m-1,1)}$.
\end{proposition} 

\begin{proof} For $m \geq 4$, consider the subspaces $M_i:=[A,[A,\ldots [L_2,L_{m-i-2}]\ldots ]/L_{m+1}$ ($0\leq i \leq m-4$) of $B_{2,m}$. They are $W_2$-submodules of $B_{2,m}$ because $W_2$ acts on $B_{2,m}$ by derivations. So the quotient space $D_{2,m}:=L_m/(L_{m+1}+M_0+\dots+M_{m-4})$ is a $W_2$-module. 
 
We claim that $D_{2,m}$ is isomorphic to $\mathcal{F}_{(m-1,1)}$ as a $W_2$-module. Take an element $[p_1,[p_2,\ldots[p_{m-1},p_m]\ldots ]$ of $D_{2,m}$. By the relations in $B_{2,3}$, we can assume that $p_m$ is either $x$ or $y$.  We notice that modulo $M_i$ we can interchange the polynomials $p_{i+1}$ and $p_{i+2}$ in the expression $[p_1,[p_2,\ldots[p_{m-1},p_m]\ldots ]$. By such permutations, we can make $p_1$ either $x$ or $y$. Similarly, using the relations in $B_{2,3}$ and permutations, we can make each of the elements $p_2, p_3,\ldots, p_{m-2}$ either $x$ or $y$. Moreover, using permutations, we can order $p_1,\ldots,p_{m-2}$ so that $p_1,\ldots,p_k=x$ and $p_{k+1},\ldots,p_{m-2}=y$ for some $0\leq k \leq {m-2}$. 

For the elements of $D_{2,m}$, we introduce the notation $c_{a,b,i,j}:=\mathrm{ad}^{a}_{x}\circ \mathrm{ad}^{b}_{y}\circ \mathrm{ad}_{x^i} (y^j)$. From the previous considerations, we obtain that $D_{2,m}[l]$ is spanned by the elements $c_{a,m-a-2,i,l-m-i+2}$, where $0\leq a \leq m-2$ and $1\leq i \leq l-m-1$. The number of these spanning elements of $D_{2,m}$ is $(m-1)(l-m-1)$.

In particular, $D_{2,m}[m]$ is spanned by the $m-1$ elements $e_i=c_{i-1,m-i-1,1,1}$, where $1\leq i \leq m-1$. We notice that in $D_{2,m}$ we have 
$$y \frac{\partial}{\partial x} e_i = \sum_{j=1}^{i-1} \mathrm{ad}^{j-1}_{x}\circ \mathrm{ad}_y (c_{i-j-1,m-i-1,1,1}) = (i-1) c_{i-2,m-i,1,1}=(i-1)e_{i-1}.$$
We notice that $e_1$ is not zero in $A$ since it has a leading monomial $xy^{m-1}$ with coefficient $(-1)^{m-2}\neq 0$. We also notice that $e_1$ has multi-degree $(1,m-1)$ in $L_m$ and $L_{m+1}[m]=M_i[m]=0$ for $0\leq i \leq m-4$. It follows that $e_1$ is not zero in the quotient space $D_{2,m}$. 

From this we derive that $e_1,e_2,\ldots,e_{m-1}$ are independent in $D_{2,m}$ because for $1\leq k \leq m-1$ if $a_k \neq 0$ we have
\begin{equation}\label{independence of e_i}
(y \frac{\partial}{\partial x})^{k-1}(\sum_{i<k}a_i e_i + a_k e_k)=a_k e_1 \neq 0
\end{equation}  
Therefore $e_1,\ldots,e_{m-1}$ form a basis of $D_{2,m}$. 

Now we show that the $W_2$-module $D_{2,m}$ is irreducible. Suppose it is not . Then it has a $W_2$-submodule $S$. Because $D_{2,m}$ starts in degree (eigenvalue of the Euler operator) $m$, $S$ has to start in degree at least $m$. We notice that the irreducible modules in the Jordan-H\"older series of $D_{2,m}$ which start in degree $m$ have the sum of their dimensions in degree $m$ equal to $m-1=\mathrm{dim} D_{2,m}[m]$. Therefore the sum of their dimensions in a degree $l>m$ will be $(m-1)(l-m+1)$. But we already showed that $\mathrm{dim} D_{2,m}[l] \leq (m-1)(l-m+1)$. Therefore all the irreducible modules in the Jordan-H\"older series of $D_{2,m}$ start in degree $m$. But the equality (\ref{independence of e_i}) shows that $D_{2,m}[m]$ belongs to a single $W_2$-submodule of $D_{2,m}$ generated by $e_1$. Therefore $D_{2,m}$ is isomorphic to an irreducible $W_2$-module, which starts in degree $m$ and has dimension $m-1$ in this degree. So this module is $\mathcal{F}_D$ where $D=(p,k)$ with $p+k=m$ and $p-k=m-2$. This is $\mathcal{F}_{(m-1,1)}$. 
\end{proof}

\section{The structures of $B_{2,3}$ and $B_{2,4}$}

In this section we find the $W_2$-module structures of $B_{2,3}$ and $B_{2,4}$. We will use characters of $W_2$-modules which are formal power series in letters $s,t$. The character of a $W_2$-module $M$ will be given by $\mathrm{char}\ M = \sum \dim M[a,b] s^a t^b$, where $M[a,b]$ denotes the subspace of elements of $M$ with weights $a,b$ of the operators $x \frac{\partial}{\partial x}, y \frac{\partial}{\partial y}$. 

First we compute the characters of the irreducible modules $\mathcal{F}_{(n,m)}$ for Young diagrams $(n,m)$.

\begin{proposition} The character of $\mathcal{F}_{(n,m)}$ is given by 
$$
\mathrm{char}\ \mathcal{F}_{(n,m)} = s^m t^m \frac{t^{n-m}+ t^{n-m-1}s + \dots + s^{n-m}}{(1-s)(1-t)}.
$$
\end{proposition}  

\begin{proof} This is true since to form an element of $\mathcal{F}_{(n,m)}[a,b]$ we have firstly to use $m$ copies of $x$ and $m$ copies of $y$ to produce the part $(dx \wedge dy)^{\otimes m}$; this accounts for the multiple $s^m t^m$ in the character formula. Next we have to choose \linebreak $0\leq i \leq n-m$ copies of $x$ and $n-m-i$ copies of $y$ to produce the symmetric part $(dx)^i \cdot (dy)^{n-m-i}$ of the tensor part of an element of $\mathcal{F}_{(n,m)}[a,b]$; this accounts for the sum $t^{n-m}+ t^{n-m-1}s + \dots + s^{n-m}$ in the numerator of the character formula. Lastly, we have to add a polynomial part to our element by multiplying it by $s^{a-m-i}$ and $t^{b-n+i}$; this is accounted for by the multiples $\frac{1}{1-s}=\sum_{l \geq 0} s^l$ and $\frac{1}{1-t}=\sum_{l \geq 0 } t^l$ in the character formula.      
\end{proof}

By multiplying $\mathrm{char}\ \mathcal{F}_{(n,m)}$ by $(1-s)(1-t)$, we obtain a polynomial with a leading monomial $s^n t^m$. Since all these polynomials for different diagrams $(n,m)$ have different leading monomials, they are independent. Therefore the characters of different $\mathcal{F}_{(n,m)}$ are linearly independent.

\begin{theorem}\label{B_{2,3} structure} 
The $W_2$-module $B_{2,3}$ is isomorphic to $\mathcal{F}_{(2,1)}$.
\end{theorem}

\begin{proof}
From the results about $B_{2,2}$ we know that $[A[A,A]A,A]\subseteq L_3$. Since $[\mathbb{C},L_2]=0$ we have that $B_{2,3}$ is a quotient of $(S(\mathbb{C}^2)/ \mathbb{C}) \otimes B_{2,2}$. By definition we have that $S(\mathbb{C}^2)$ is isomorphic to $\mathcal{F}_{(0,0)}$. By the results of \cite{FS} we also have that $B_{2,2}$ is isomorphic to $\mathcal{F}_{(1,1)}$. So $B_{2,3}$ is a quotient of $(\mathcal{F}_{(0,0)} / \mathbb{C}) \otimes \mathcal{F}_{(1,1)}$. Therefore the irreducible modules in the Jordan-H\"older series of $B_{2,3}$ will be found among the irreducible modules contained in the module $(\mathcal{F}_{(0,0)} / \mathbb{C}) \otimes \mathcal{F}_{(1,1)}$. To find them, we compute the character of the last module:
\begin{eqnarray*}
&&\mathrm{char}\ (\mathcal{F}_{(0,0)}/ \mathbb{C}) \otimes \mathcal{F}_{(1,1)}\\ 
&=& (\frac{1}{(1-s)(1-t)}-1) \frac{st}{(1-s)(1-t)}\\
&=&\sum_{k \geq 0} st \frac{(s^k+s^{k-1}t+\dots+t^k)}{(1-s)(1-t)} - \mathrm{char}\ \mathcal{F}_{(1,1)}\\
&=& \sum_{p \geq 1} \mathrm{char}\ \mathcal{F}_{(p,1)} - \mathrm{char}\ \mathcal{F}_{(1,1)} = \sum_{p \geq 2} \mathrm{char}\ \mathcal{F}_{(p,1)}.
\end{eqnarray*}
But we know from Theorem \ref{multiplicities} that the only copy of $\mathcal{F}_{(p,1)}$ in $B_{2,3}$ is $\mathcal{F}_{(2,1)}$. Therefore $B_{2,3}$ is isomorphic to $\mathcal{F}_{(2,1)}$.   
  
\end{proof}

\begin{theorem}\label{B_{2,4} structure} The $W_2$-module $B_{2,4}$ has in its Jordan-H\"older series only two irreducible $W_2$-modules, $\mathcal{F}_{(3,1)}$ and $\mathcal{F}_{(3,2)}$ and each with multiplicity $1$.
\end{theorem}

\begin{proof}
From the results about $B_{2,2}$, we know that $[A[A,A]A,A]\subseteq L_3$.  Since $[\mathbb{C},L_3]=0$ we have that $B_{2,4}$ is a quotient of $(S(\mathbb{C}^2)/ \mathbb{C}) \otimes B_{2,3}$. 

Since $S(\mathbb{C}^2)$ is isomorphic to $\mathcal{F}_{(0,0)}$ and from Theorem \ref{B_{2,3} structure}, we know that $B_{2,4}$ is a quotient of $(\mathcal{F}_{(0,0)} / \mathbb{C}) \otimes \mathcal{F}_{(2,1)}$. Therefore the irreducible modules in the Jordan-H\"older series of $B_{2,4}$ will be found among the irreducible modules contained in $(\mathcal{F}_{(0,0)} / \mathbb{C}) \otimes \mathcal{F}_{(2,1)}$. By a similar computation to the one in the proof of Theorem \ref{B_{2,3} structure}, we have
$$
\mathrm{char}\ (\mathcal{F}_{(0,0)}/ \mathbb{C}) \otimes \mathcal{F}_{(2,1)} =\sum_{p \geq 3} \mathrm{char}\ \mathcal{F}_{(p,1)} + \sum_{p \geq 2} \mathrm{char}\ \mathcal{F}_{(p,2)}.
$$

But we know from Theorem \ref{multiplicities} that the only copy of $\mathcal{F}_{(p,1)}$ in $B_{2,4}$ is $\mathcal{F}_{(3,1)}$ and the only copy of $\mathcal{F}_{(p,2)}$ is $\mathcal{F}_{(3,2)}$. Therefore, the Jordan-H\"older series of the module $B_{2,4}$ contains exactly two irreducible $W_2$-modules $\mathcal{F}_{(3,1)}$ and $\mathcal{F}_{(3,2)}$.   
  
\end{proof}

Now we show that $B_{2,4}$ is not a direct sum of the modules $\mathcal{F}_{(3,1)}$ and $\mathcal{F}_{(3,2)}$ in its Jordan-H\"older series.

\begin{proposition} The $W_2$-module $B_{2,4}$ is isomorphic to a nontrivial extension of $\mathcal{F}_{(3,2)}$ by $\mathcal{F}_{(3,1)}$. 
\end{proposition}

\begin{proof} 
For a $W_2$-module $M$, we denote by $M[k]$ the weight space of $M$ for weight $k$ of the Euler vector field in two variables.  Notice that $B_{2,4}$ has a $W_2$-submodule $C_{2,4} = [[A_2,A_2],[A_2,A_2]]/L_5$. The lowest weight of $C_{2,4}$ is $5$ and the lowest weight vectors are $a[[x^2,y],[x,y]]+b[[x,y^2],[x,y]]$. Therefore $C_{2,4}$ is isomorphic to $\mathcal{F}_{(3,2)}$. Since we have $\dim \mathcal{F}_{(3,2)}[4] = 2$, it follows that $[[x^2,y],[x,y]]$ and $[[x,y^2],[x,y]]$ form a basis of $C_{2,4}[4]$. From Theorem \ref{B_{2,4} structure} it follows that the $W_2$-module $B_{2,4}/C_{2,4}$ is isomorphic to $\mathcal{F}_{(3,1)}$. So we have an exact sequence of $W_2$-modules
$$
0 \rightarrow \mathcal{F}_{(3,2)} \rightarrow B_{2,4} \rightarrow \mathcal{F}_{(3,1)} \rightarrow 0.
$$
 
We will now show that this sequence does not split. 
Since the diagram $(3,2)$ has $5$ cells, $\dim \mathcal{F}_{(3,2)}[4]=0$. Then if we had $B_{2,4} \cong \mathcal{F}_{(3,1)} \oplus \mathcal{F}_{(3,2)}$, the entire space $B_{2,4}[4]$ would belong to the copy of $\mathcal{F}_{(3,1)}$ in $B_{2,4}$ which we denote by $F$. Notice that $[x,[x,[x,y]]]$, $[x,[y,[x,y]]]$, $[y,[y,[x,y]]]$ are in $B_{2,4}[4]$, so 
$$s=-3y^2 \frac{\partial}{\partial y} [x,[x,[x,y]]] -x^2 \frac{\partial}{\partial y} [y,[y,[x,y]]] + 2xy \frac{\partial}{\partial x} [x,[x,[x,y]]]$$ is in $F$. 

By using Jacobi identity and relations in $B_{2,3}$, we have:
\begin{eqnarray*}
-3y^2 \frac{\partial}{\partial y} [x,[x,[x,y]]]&=&-3 [x,[x,[x,y^2]]],\\
-x^2 \frac{\partial}{\partial y} [y,[y,[x,y]]]&=& [[x,y],[x^2,y]] -2[y,[x^2,[x,y]]],\\
2xy \frac{\partial}{\partial x} [x,[x,[x,y]]]& =& 
 [[x,y],[x^2,y]] + 2[x,[y,[x^2,y]]] + 3[x,[x,[x,y^2]]]. 
\end{eqnarray*}

Adding them up, we obtain
$s=4[[x,y],[x^2,y]]$ which is a nonzero element in $B_{2,4}$.
Since $s$ belongs to $F \cap C_{2,4}$, we have that $F \cap C_{2,4} \neq 0$ which contradicts our assumption that $B_{(2,4)} = F \oplus C_{2,4}$. So as a $W_2$-module, $B_{2,4}$ is isomorphic to a nontrivial extension of $\mathcal{F}_{(3,2)}$ by $\mathcal{F}_{(3,1)}$.  
\end{proof}

To completely characterize $B_{2,4}$ as a $W_2$-module, we prove

\begin{proposition} All the nontrivial extensions of $\mathcal{F}_{(3,2)}$ by $\mathcal{F}_{(3,1)}$ are isomorphic.
\end{proposition}

\begin{proof}
Firstly we construct such a nontrivial extension abstractly. We
have the Lie algebra $W_n$ of polynomial vector fields on
$V^*$, where $V=\mathbb{C}^n$. We denote by $W_n^0$ the subalgebra of $W_n$ of
vector fields vanishing at the origin. For every Young diagram
$D$, we have a corresponding representation $F_D$ of
$\mathfrak{gl}(n,\C)$, and a corresponding representation of
$W_n^0$ in which linear vector fields $\sum a_{ij} x_i
\frac{\partial}{\partial x_j}$ act as matrices $(a_{ij})$ and
higher-order vector fields act by zero. Suppose that $D, E$ are
two Young diagrams such that if we align their left upper corners
the set-theoretic difference $E-D$ is equal to one box (an
example of such a pair of diagrams is $E=(3,2)$, $D=(3,1)$). It
is known that in this case there exists a nonzero homomorphism
$F_D \otimes V \rightarrow F_E$, which is unique up to scaling. 

We construct a representation $Y$ of $W_n^0$ as follows. As a
vector space $Y:=F_D \oplus F_E$. Linear vector fields which
correspond to $\mathfrak{gl}(n,\C)$ act on $Y$ as in the direct
sum of the representations $F_D, F_E$ of
$\mathfrak{gl}(n,\C)$. Cubic and higher vector fields act by
zero. It remains to describe how quadratic vector fields
act. They form a space $S^2V\otimes V^*$, which has a unique
invariant projection to $V$. So we can define an action of
$S^2V\otimes V^*$ on $F_D\oplus F_E$ by using this projection and the map $F_D \otimes V \rightarrow F_E$ (this action will map the subspace $F_D$ to $F_E$ and the subspace $F_E$ to $0$).

Now we define the representation $\mathcal{F}_Y:=\mathrm{Hom}_{U(W_n^0)}(U(W_n),Y)$. Then we have an exact sequence   
$$0\to \mathcal{F}_E \rightarrow \mathcal{F}_Y \rightarrow
\mathcal{F}_D\to 0.$$
From now on, let us fix the Young diagrams $D=(3,1)$, $E=(3,2)$ and the corresponding representations $Y,\mathcal{F}_Y$ of $W_2^0$ and $W_2$. 

Now we prove that any $W_2$-module $M$ for which there is a short exact sequence
$$0\to \mathcal{F}_{(3,2)} \rightarrow M \rightarrow
\mathcal{F}_{(3,1)}\to 0$$
which does not split is isomorphic to $\mathcal{F}_Y$. Suppose we
have such a module $M$. We have
$M[4]\cong F_{(3,1)}$ and $M[5]\cong F_{(3,2)} \oplus
F_{(3,1)} \otimes V$, which is isomorphic to
$F_{(3,2)}(1) \oplus (F_{(3,2)}(2) \oplus F_{(4,1)})$
(the 1 and 2 in parentheses denote the first and the second copy). So
as $\mathfrak{sl}(2,\C)$-modules, $M[4]\cong V_2$ and
$M[5]\cong V_1(1) \oplus (V_1(2) \oplus V_3)$ where the subscripts
denote the highest weights. Now we have the degree $1$ part $W[1]$
(quadratic vector fields) of $W:=W_2$ acting from $M[5]^*$ to
$M[4]^*$. As an $\mathfrak{sl}(2,\C)$-module, we have a
decomposition $W[1]=V_1\oplus V_3$. Let us pick a nonzero
element $f$ in $(V_1(1)\oplus V_1(2))^*\subset M[5]^*$ of weight $1$ which is killed by the
lowest vector (of weight $-3$) of $V_3\subset W[1]$. This is a
scalar linear equation, so $f$ exists (and is unique up to a
scalar since the above equation is nontrivial). It generates a
copy of $V_1$ inside $M[5]^*$, which we call $N$. Moreover, since
the extension is nontrivial, $W[1]$ acts nontrivially on
$N$. Thus, $N^\perp\oplus M[\ge 6]\subset M$ is a
$W_2^0$-submodule, and the quotient module $M/(N^\perp\oplus
M[\ge 6])=N^*\oplus M[4]$ is isomorphic to $Y$. 

Therefore we have a natural $W_2$-homomorphism $M\to {\mathcal
F}_Y$, which is an isomorphism in degrees $4$ and $5$. Hence it
is an isomorphism (as there are only $2$ terms in the
Jordan-H\"older series of $M$). 
The proposition is proved. 
\end{proof}

{\bf Acknowledgments}

We thank Prof. Pavel Etingof for giving this problem and for many
useful discussions. We also thank Prof. Ju-Lee Kim for many
useful discussions and Jennifer Balakrishnan for much help with
the software ``Magma''. The work of J.K. and X.M. was done 
within the framework of SPUR (Summer Program of Undergraduate
Research) at the Mathematics Department of MIT in July 2007. 

\section{Appendix: $B_2(A)$ for a general associative algebra $A$}
\centerline{by Pavel Etingof} 

\vskip .05in

The goal of this appendix is to generalize some 
of the results of Feigin and Shoikhet \cite{FS} to the case of 
any associative algebra. 

\subsection{The algebra $R(A)$}
Let $A$ be an associative algebra over $\C$. 
Let $D(A)=A\oplus A$, regarded as a supervector space, where the first copy of 
$A$ is even and the second one is odd. For $a\in A$, let us denote the 
elements $(a,0),(0,a)$ of $D(A)$ by $x_a,\xi_a$, respectively. 
 
Define the supercommutative algebra $R(A)$ to be the quotient of 
the symmetric algebra $SD(A)$ by the relations 
$$
x_ax_b-x_{ab}+\xi_a\xi_b=0
$$
and 
$$
x_a\xi_b+\xi_ax_b-\xi_{ab}=0.
$$
This is a DG algebra, with $dx_a=\xi_a,d\xi_a=0$. 

It is clear that the quotient of $R(A)$ by the ideal $I$ generated by 
the odd elements is 
$A_{\rm ab}$, the abelianization of $A$. Thus $R(A)$ is a certain 
super-extension of the abelianization of $A$. 
More precisely, let $\Omega(A_{\rm ab})$ be the DG algebra of 
K\"ahler differential forms
for the abelianization $A_{\rm ab}$ of $A$.  
It is defined by the same generators as $R(A)$ but 
with defining relations 
$$
x_ax_b-x_{ab}=0
$$
and 
$$
x_a\xi_b+\xi_ax_b-\xi_{ab}=0
$$
with $dx_a=\xi_a$, $d\xi_a=0$. 
Thus, denoting by ${\rm gr}R(A)$ the associated graded algebra of $R(A)$ 
under the filtration by powers of $I$, we obtain that 
there is a natural surjective homomorphism $\eta: \Omega(A_{\rm ab})\to 
{\rm gr}R(A)$. It is not always an isomorphism.  

\begin{definition}\label{smoo}
We will say that $A$ is pseudosmooth if 
$A_{\rm ab}$ is a regular finitely generated algebra 
(i.e. ${\rm Spec}(A_{\rm ab})$ is a smooth affine algebraic variety $X$), and 
$\eta$ is an isomorphism.
\end{definition} 

\begin{proposition}\label{abb}
$A$ is pseudosmooth if and only if $R(A)$ is isomorphic, as a DG
algebra, to the algebra $\Omega(X)$ of regular differential forms
on a smooth affine algebraic variety $X$. 
\end{proposition}

\begin{proof}
Suppose that $R(A)=\Omega(X)$. Then
$A_{\rm ab}=\Omega(X)/(d\Omega(X))={\mathcal O}_X$, 
and $\eta$ is clearly an isomorphism. Conversely, if $A_{\rm
ab}={\mathcal O}_X$ for smooth $X$ and $\eta$ is an isomorphism then 
the projection $R(A)\to A_{\rm ab}$ splits, and this splitting 
uniquely extends to an isomorphism of DG algebras $\Omega(X)\to
R(A)$. 
\end{proof}   

\subsection{The Fedosov products}

For any DG algebra $S$ introduce the Fedosov product on $S$ 
by 
$$
f*g=f\cdot g+(-1)^{|f|}df\cdot dg,
$$
and the inverse Fedosov product by 
$$
f\circ g=f\cdot g-(-1)^{|f|}df\cdot dg,
$$
and let $S_*,S_\circ$ be the algebra $S$ equipped with the Fedosov product,
respectively the inverse Fedosov product. 

Obviously, the operations of passing to the Fedosov and inverse Fedosov product
in a differential algebra are inverse to each other, hence the terminology. 

\subsection{The universal property} 

It turns out that the algebra $R(A)$ has the following universal property. 

\begin{proposition}\label{unprop}
For any supercommutative DG algebra $S$, one has 
a natural isomorphism ${\rm Hom}_{DG}(R(A),S)\to {\rm Hom}(A,S_{*0})$, where 
$S_{*0}$ is the even part of $S_{*}$. 
\end{proposition} 

\begin{proof}
It is clear that any homomorphism
$f: R(A)\to S$ is determined by the elements $y_a=f(x_a)$, and 
the elements $y_a$ define a homomorphism if and only if they satisfy the 
equations $y_a*y_b-y_{ab}=0$. This implies the statement. 
\end{proof} 

\subsection{Relation with noncommutative differential forms}

In fact, the algebra $R(A)$ can be obtained from noncommutative 
differential forms on $A$ (\cite{CQ}). 
Namely, let $\Omega_{\rm nc}(A)=A\otimes T(\bar A)$ denote the DG algebra 
of noncommutative differential forms on $A$ (here $\bar A=A/\Bbb C$);
it is the span of formal expressions $a_0da_1\cdots da_n$. 

\begin{proposition}\label{ncd}
The algebra $R(A)$ is naturally isomorphic to 
the abelianization (in the supersense) of the DG algebra 
$\Omega_{\rm nc}(A)_\circ$.  
\end{proposition}

\begin{proof}  
It suffices to show that if 
$S$ is a supercommutative DG algebra,\linebreak then 
${\rm Hom}_{DG}(R(A),S)={\rm Hom}_{DG}(\Omega_{\rm nc}(A)_\circ,S)$. 

But ${\rm Hom}_{DG}(\Omega_{\rm nc}(A)_\circ,S)=
{\rm Hom}_{DG}(\Omega_{\rm nc}(A),S_*)={\rm Hom}_{DG}(A,S_{*0})$, and 
the result follows from the universal property of $R(A)$.   
\end{proof}  

\subsection{Description of $R(A)$ using a presentation of $A$}

Let $V$ be a vector space. Then $SV\otimes \wedge V$ is 
naturally a differential algebra (the De Rham complex of $V^*$). 
Suppose that $A=TV/(L)$, where $L\subset TV$ is a set of relations. 

Let $g: TV\to (SV\otimes \wedge V)_{*0}$ 
be the homomorphism defined by the condition that $g(v)=v\in SV$ 
for $v\in V$. 

\begin{proposition}\label{rela} We have
$R(A)=(SV\otimes \wedge V)/(g(L)\cup dg(L))$. 
\end{proposition}

In particular, we see that $R(TV)=SV\otimes \wedge V$. 

\begin{proof}
We have 
$$
{\rm Hom}(A,S_{*0})=
\lbrace{f\in {\rm Hom}_{DG}(SV\otimes \wedge V,S): f(g(L))=0\rbrace},
$$ 
which implies the desired statement by Proposition \ref{unprop}.
\end{proof}

\subsection{The quotient of $A$ by triple commutators}

\begin{proposition}\label{tripcom}
We have a natural isomorphism of algebras 
$$
\phi: A/A[[A,A],A]A\to R(A)_{*0}.
$$  
\end{proposition}

\begin{proof}
We have a natural homomorphism $\phi$ given by $\phi(a)=x_a$.
Let us show that it is an isomorphism. As shown in \cite{FS}, $\phi$ 
is an isomorphism for $A=TV$. On the other hand, $A/A[[A,A],A]A$ 
is the quotient of $TV/TV[[TV,TV],TV]TV$ by the additional relations $L$. 
Thus, it suffices to show that $R(A)_{*0}$ is obtained from 
$(SV\otimes \wedge V)_{*0}$ by imposing additional relations
$g(L)$. These relations clearly hold, so we need to show that 
there is no others. 

Thus, by Proposition \ref{rela}, we need to show that in the algebra 
$(SV\otimes \wedge V)_{*0}/(g(L))$, we have 
$a\cdot g(b)=0$ and $c\cdot dg(b)=0$ for all $b\in L$, 
$a\in (SV\otimes \wedge V)_0$ and $c\in (SV\otimes \wedge V)_1$.
 
The first equality follows since 
$a\cdot g(b)=\frac{1}{2}(a*g(b)+g(b)*a)$.  
To prove the second equality, note that since $c$ is odd, 
we have $c=\sum c_j\cdot dv_j$, hence\linebreak $c\cdot dg(b)=\sum c_j\cdot 
dv_j\cdot dg(b)$, and $dv\cdot dg(b)=\frac{1}{2}(v*g(b)-g(b)*v)=0$. 
\end{proof}

\begin{proposition}\label{imd}
The map $\phi$ of Proposition \ref{tripcom} 
maps $[A,A]$ onto the image of $d$ in $R(A)_{*0}$.
\end{proposition}

\begin{proof}
It is shown in \cite{FS} that if $A=F$ is a free algebra 
then the statement holds. This implies that it holds for any
associative algebra. 
\end{proof}

Let ${\rm gr}(A)$ be the associated graded Lie algebra of $A$ 
with respect to its lower central series filtration.
Let $Z(A)=A[[A,A],A]A/([A,A]\cap A[[A,A],A]A)$. 
Thus, $Z(A)\subset B_1(A)$.

\begin{proposition}\label{cente} 
(i) $Z(A)$ is central in the Lie algebra ${\rm gr}(A)$.  

(ii) The space $B_1(A)/Z(A)$ is isomorphic, via $\phi$, to 
$R(A)_0/R(A)_0^{\rm exact}$. 
\end{proposition}

\begin{proof}
Part (i) follows from Lemma 2.2.1 of \cite{FS}
(this lemma is proved in \cite{FS} for the free algebra 
but applies without changes to any associative algebra).
Part (ii) follows from Proposition \ref{imd}.  
\end{proof} 

\subsection{The first cyclic homology}

Let $A$ be an associative algebra, and 
$W(A)$ be the subspace of $\wedge^2A$ 
spanned by the elements 
$$
ab\wedge c+bc\wedge a+ca\wedge b.
$$
We have a natural map 
$[,]: \wedge^2A/W(A)\to [A,A]$ given by $a\wedge b\to [a,b]$. 
Recall \cite{Lo} that the first cyclic homology $HC_1(A)\subset
\wedge^2A/W(A)$ is the kernel of this map. 

Define the map $\zeta: \wedge^2A/W(A)\to R(A)_1/R(A)_1^{\rm exact}$ by
the formula 
$$
\zeta(a\wedge b)=d\phi(a)\cdot \phi(b).
$$
It is easy to see that this map is well defined. 
Moreover, if $u\in HC_1(A)$ then $\zeta(u)$ closed. 
Thus, we obtain a map $\zeta: HC_1(A)\to H^{\rm odd}(R(A))$. 
Denote by $Y(A)$ the image of this map. 

\subsection{Pseudoregular DG algebras}

Let $S$ be a commutative DG algebra. 
Let $S'=S/S^{\rm exact}$. Define the linear map 
$\theta: \wedge^2S_0\to S_1'$
by the formula $\theta(a,b)=da\cdot b$. This is skew-symmetric
because 
$da\cdot b+db\cdot a=d(ab)$. It is clear that the kernel 
$\ker\theta$ contains the elements 
$$
\kappa(a,b,c):=ab\wedge c+bc\wedge a+ca\wedge b,
$$
where $a,b,c\in S_0$, and the elements $a\wedge b$ 
where $a$ is exact. Denote the span of these two types of
elements by $E$.  

Let us say that $S_0$ is {\it pseudoregular} if 
$S_1=S_0dS_0$ (implying that $\theta$ is surjective), and 
$\ker\theta=E$. 

\subsection{Pseudoregularity of the De Rham DG algebra of a smooth variety}

Let $X$ be a smooth affine algebraic variety over $\C$. 
Denote by ${\mathcal O}_X$ the algebra of regular functions on
$X$, and by 
$\Omega(X)$ the DG algebra of regular differential forms on $X$. 

\begin{theorem}\label{pseudore}
The algebra $S:=\Omega(X)$ is pseudoregular. 
\end{theorem}

\begin{remark}
This was proved in \cite{FS} in the special case when $X$ is the
affine 
space $\C^n$. 
\end{remark}

\begin{proof} It is obvious that $S_1=S_0dS_0$. 
We need to show that $\theta$ identifies 
$\wedge^2S_0/E$
with $S_1'$. 
To do so, write $\kappa(a,b,c)$ in the form 
$$
\kappa(a,b,c)=ab\wedge c-a\wedge bc-b\wedge ca.
$$
From this we see that modulo the span of $E$, 
any element of $\wedge^2S_0$
can be reduced to an element of 
${\mathcal O}_X\otimes S_0$ (where ${\mathcal O}_X$ is viewed as
the 
subspace of 0-forms in the space $S_0$ of even forms). 

Furthermore, by modding out by $\kappa(a,b,c)$ we factor out a
subspace of \linebreak 
${\mathcal O}_X\otimes S_0$ which is spanned by $ab\otimes
g-a\otimes bg-
b\otimes ga$, $a,b\in {\mathcal O}_X$, $g\in S_0$. The
corresponding 
quotient space is the Hochschild homology $HH_1({\mathcal
O_X},S_0)$. 
Since $S_0$ is a projective module over ${\mathcal O}_X$ (as $X$
is smooth),
we have $HH_1({\mathcal O_X},S_0)=HH_1({\mathcal O}_X,{\mathcal 
O}_X)\otimes_{{\mathcal O}_X}S_0$, which by the
Hochschild-Kostant-Rosenberg
theorem (\cite{Lo}) equals $\Omega^1(X)\otimes_{{\mathcal O}_X}
S_0$. 
In fact, the relevant projection ${\mathcal O}_X\otimes_\C S_0\to
\Omega^1(X)\otimes_{{\mathcal O}_X} S_0$ is simply given by the
formula
$a\otimes g\to da\otimes g$. 

Further, for any $a,b,c\in {\mathcal O}_X$, $f\in S_0$ we have,
modulo $E$: 
$$
a\wedge db\cdot dc\cdot f=a\cdot db\cdot dc\wedge f=-b\cdot
da\cdot dc\wedge 
f=-b\wedge da\cdot dc\cdot f,
$$
which proves that in fact, modulo $E$, the
space 
$\Omega^1(X)\otimes_{{\mathcal O}_X} S_0$ gets projected onto its
quotient 
space $S_1$, and the resulting projection map 
${\mathcal O}_X\otimes_\C S_0\to S_1$ is given by $a\otimes g\to
da\cdot g$. 
Moreover, it is clear that we project further down to
$S_1'=S_1/S_1^{\rm 
exact}$, because the space of exact elements of $S_1$ is spanned
by 
elements of the form $da\wedge f$, where $f$ is an exact element
of $S_0$ and $a\in {\mathcal O}_X$, 
and such an element is the image of 
$a\wedge f$, which belongs to $E$. The theorem is proved. 
\end{proof}

\subsection{The structure of $B_2(A)$ for pseudosmooth algebras}

The main result of the appendix is the following theorem. 

\begin{theorem}\label{pseudosm}
Let $A$ be a pseudosmooth algebra.
Then
 
(i) $B_2(A)$ is naturally isomorphic 
to $R(A)_1'/Y(A)$; in particular, if $R(A)$ has no odd
cohomology, then $B_2(A)=R(A)_0^{\rm exact}$. 

(ii) $([A,A]\cap A[[A,A],A]A)/[[A,A],A]$ 
is naturally isomorphic to $H^{\rm odd}(R(A))/Y(A)$. 

(iii) In terms of the identification of (i) and Proposition \ref{imd}(ii),
the bracket map $\wedge^2(B_1(A)/Z(A))\to B_2(A)$
is given by the formula $a\wedge b\to da\cdot b$. 
\end{theorem}

\begin{proof}
According to \cite{FS}, proof of Lemma 1.2, we have an exact
sequence
$$
\cdots\to HC_1(A)\to \wedge^2(A/[A,A])/(ab\wedge c+bc\wedge
a+ca\wedge b)\to [A,A]/[[A,A],A]\to 0.
$$
By Proposition \ref{cente}, this implies that we have an exact sequence 
$$
\cdots\to HC_1(A)\to \wedge^2(R(A)'_0)/(ab\wedge c+bc\wedge
a+ca\wedge b)\to [A,A]/[[A,A],A]\to 0.
$$
Since $A$ is pseudosmooth, by Proposition \ref{abb} 
the middle term has the form 
$$
\wedge^2\Omega_{\rm even}'(X)/(a\wedge bc+b\wedge
ca+c\wedge ab),
$$ 
where $X$ is the spectrum of $A_{\rm ab}$. 
By Theorem \ref{pseudore},
this equals $\Omega_{\rm odd}(X)/\Omega_{\rm odd}^{\rm
exact}(X)$. Clearly, the space $HC_1(A)$ maps onto 
$Y(A)\subset \Omega_{\rm odd}(X)/\Omega_{\rm odd}^{\rm
exact}(X)$. This implies the first and third statements. 
The second statement follows from the first one and 
Proposition \ref{imd}. 
\end{proof}

\begin{remark}
In the special case when $A$ is a free algebra, 
Theorem \ref{pseudosm} is proved in \cite{FS}. 
In this case, one has $[A,A]\cap A[[A,A],A]A=[[A,A],A]$.
However, in general this equality does not have to hold. 
For example, let $A$ be the algebra generated by two elements
$x,y$ with the only relation $xy=1$. Then it is easy to show that
$HC_1(A)=0$ (see e.g. \cite{EG}, Section 5.4), and $R(A)=
\Omega(X)$, where $X$
is the curve defined by the equation $xy=1$ in the plane (i.e
$X=\Bbb C^*$). This algebra is commutative (even with the
$*$-product), since $X$ is 1-dimensional. Thus, 
$A/A[[A,A],A]A$ is commutative, and hence $[A,A]\subset
A[[A,A],A]A$. However, it follows from Theorem 
\ref{pseudosm} that the space $[A,A]/[[A,A],A]$ is 1-dimensional.
In fact, one may check that it is spanned by the element $[x,y]$. 
\end{remark} 

\subsection{A sufficient condition of pseudosmoothness}

\begin{proposition}\label{pseudos}
Let $f_1,\ldots,f_m\in A_n$ be a set of elements,  
such that their images $\bar f_1,\ldots,\bar f_m$ in 
$\Bbb C[x_1,\ldots,x_n]$ form a regular sequence defining a 
smooth complete intersection $X$ in $\Bbb C^n$ (of codimension
$m$). Then the algebra $A:=A_n/(f_1,\ldots,f_m)$ is pseudosmooth, 
and $R(A)$ is isomorphic to $\Omega(X)$.    
\end{proposition}

\begin{proof}
We have $A_{\rm ab}={\mathcal O}_X$, and because 
of the complete intersection condition, $\eta$ is an
isomorphism. Thus $A$ is pseudosmooth, and by 
Proposition \ref{abb}, $R(A)$ isomorphic to
$\Omega(X)$. 
\end{proof} 

\subsection{Examples} 

The above results allow one to compute $B_2(A)$ for specific algebras $A$. 

\begin{proposition}\label{kah} Suppose that $L\subset SV\subset TV$. Then 
$g(L)=L\subset SV$, and hence $R(A)$ is naturally isomorphic to 
the algebra of K\"ahler differential forms $\Omega(A_{\rm ab})$. 
In particular, if $A_{\rm ab}$ is regular then $A$ is
pseudosmooth. 
\end{proposition}

\begin{proof} Obvious. 
\end{proof}

\begin{example}\label{sl2} 
Let $A$ be the free algebra in three generators $x,y,z$ modulo the relation 
$x^2+y^2+z^2=1$ (noncommutative 2-sphere). Let us compute the 
space $B_2(A)$ as a representation of $SO(3)$ acting on this algebra. 
From Proposition \ref{kah} we find that $A$ is pseudosmooth, and 
$R(A)$ is the algebra of polynomial 
differential forms on the usual commutative quadric $Q$. In this
case, we have no odd cohomology, so by Theorem \ref{pseudosm}, 
$B_2(A)$ is the space of exact 2-forms. The space of exact 2-forms is a 
subspace of codimension 1 in the space of all 2-forms, since 
$H^2(Q)$ is 1-dimensional. The space of all 
2-forms is isomorphic to the space of functions as an 
$SO(3)$-module, since there 
is an invariant symplectic form on the quadric (the area form).  
Now, we have 
$$
{\rm Fun}(Q)=V_0\oplus V_2\oplus V_4\oplus\cdots,
$$
where $V_{2i}$ is the $(2i+1)$-dimensional representation of $SO(3)$. 
Thus, 
$$
B_2(A)=V_2\oplus V_4\oplus\cdots
$$
\end{example} 

Let us now consider more general examples.  
As before, assume that $L\subset SV$,
and suppose that $A=TV/(L)$ is a pseudosmooth algebra
(i.e., $A_{\rm ab}$ is regular),
such that $R(A)$ has no odd cohomology. 
Suppose further that $L$ is fixed by a reductive subgroup $G\subset
GL(V)$, such that $R(A)$ is a direct sum of irreducible
representations of $G$ with finite multiplicities. In this case, 
one can define the character-valued Hilbert series
$F(z),E(z),H(z)$ of the graded representations
$R(A)$, $R(A)^{\rm exact}$, and the cohomology $H(A)$ of $R(A)$. 
Then we have the equations 
$$
z(F-E-H)=E, 
$$
which implies that 
\begin{equation}\label{form-e}
E=\frac{z(F-H)}{1+z}.
\end{equation}
This formula is useful because often $F$ and $H$ are known
explicitly. 

\begin{example} Let $\g$ be a simple Lie algebra with root system
$R$ and Weyl group $W$, and let $G$ the corresponding simply connected group. 
Let $r$ be the rank of $G$, 
$p_1,\ldots,p_r$ be homogeneous generators of the ring $(S\g)^G$,
and $d_i=\deg(p_i)$. Let $b_i$ be generic complex numbers, 
and let $A(\g,b)$ be the
quotient of the tensor algebra $T\g$ by the relations $p_i=b_i$.
Note that the algebra from Example \ref{sl2} is the special case
of $A(\g,b)$ for $\g={\mathfrak{sl}}(2)$. 

Let us calculate the 
decomposition of the space $B_2(A)$ into irreducible
representations of $G$. We have $B_2(A)=\oplus_{V\in {\rm
Irr}(G)} N_V\otimes V$, where $N_V={\rm Hom}_G(V,B_2(A))$. 

By formula (\ref{form-e}), we have  
$$
\dim N_V=\frac{1}{2}(E_V(1)+E_V(-1)),
$$
with
$$
E_V(z)=\frac{z}{1+z}(F_V(z)-H_V(z)),
$$
where $F_V$ and $H_V$ are contributions of $V$ into $F$ and $H$,
respectively. It remains to find $F_V(z)$ and $H_V(z)$. 

By Proposition \ref{kah}, we find that $R(A)$ is 
the algebra of polynomial differential forms 
on $G/H$, where $H$ is a maximal torus in $G$. 
Thus we have 
$H_V(z)=0$ unless $V=\C$, 
$$
H_\C(z)=\prod_{i=1}^r \frac{z^{2d_i}-1}{z^2-1}
$$
is the Poincar\'e polynomial of $G/H$, and 
$$
F_V(z)=\sum_{j\ge 0}z^j\dim {\rm Hom}_H(V,\wedge^j(\g/\h)),
$$
where $\h={\rm Lie}H$.
More explicitly, 
$$
F_V(z)=C.T.(\chi_{V^*}\cdot \prod_{\alpha\in R}(1+ze^\alpha)),
$$
where $\chi_{V^*}$ is the character of $V^*$, and C.T. means the
constant term.  

In the case $\g={\mathfrak{sl}}(2)$, this recovers the answer from 
Example \ref{sl2}. 
 
\begin{corollary}
Let $\nu(R)$ be the number of subsets of $R$ with zero sum. 
Then 
$$
\dim B_2(A)^G=\frac{1}{4}(\nu(R)-|W|).
$$ 
\end{corollary} 

\begin{proof} It is easy to show that $F_\C(-1)=H_\C(-1)=|W|$, and 
$F_\C'(-1)=H_\C'(-1)=-|R||W|/2$, thus $E_\C(-1)=0$. 
So $\dim B_2(A)^G=\frac{1}{2}E_\C(1)=\frac{1}{4}(F_\C(1)-H_\C(1))$. 
But we have $H_\C(1)=|W|$, and $F_\C(1)=\nu(R)$. 
The corollary follows. 
\end{proof}

For example, $\dim B_2(A)^G$ is $0$ for $\g={\mathfrak{sl}}(2)$, 
$1$ for $\g={\mathfrak{sl}}(3)$, and $32$ for $\g={\mathfrak{sl}}(4)$. 
\end{example} 

\begin{example} Let $P\in \C\langle x,y\rangle$ 
be a noncommutative polynomial in two
variables $x,y$, and $\bar P$ be the abelianization of $P$ 
(i.e., the image of $P$ in the polynomial algebra $\C[x,y]$). 
Denote by $A_P$ the algebra $\C\langle x,y\rangle/(P)$. 
Assume that the curve $X_{\bar P}$ given by the equation $\bar P(x,y)=0$
is smooth. Then $A=A_P$ is pseudosmooth, and thus Theorem
\ref{pseudosm} applies to $A$. Moreover, since $X_{\bar P}$ is a curve, 
the algebra $A/A[[A,A],A]A$ is commutative, and hence 
$[A,A]\subset A[[A,A],A]A$. Thus $B_2(A)=([A,A]\cap
A[[A,A],A]A)/[[A,A],A]=H^1(X_{\bar P})/Y(A)$. 

The space $Y(A)$ actually depends on $P$, not only on $\bar P$. 
For example, assume that the leading term of $P$ is generic.
In this case, by the results \cite{EG}, $HC_1(A)=0$, 
and hence $B_2(A)=H^1(X_{\bar P})$. The same holds if 
the leading term of $P$ is, say $x^py^q$. Thus, for example, if 
$P=x^2y-1$ then $B_2(A)=H^1(\C^*)=\C$. On the other hand, if 
$P=xyx-1$ then in $A$ we have $xy=xyxyx=yx$, so
$A=\C[x,y]/(yx^2=1)$, and $B_2(A)=0$ (thus, $Y(A)$ is
1-dimensional in this case). 

Let us do two concrete examples. 

1. $P$ is a generic polynomial of degree $d$. In this case 
the curve $X=X_{\bar P}$ has genus $(d-1)(d-2)/2$ and $d$ points
at infinity. So its Euler characteristic is 
$\chi=2-(d-1)(d-2)-d=-d(d-2)$, and hence 
$\dim B_2(A)=\dim H^1(X)=(d-1)^2$.

2. Let $P(x,y)=Q(x)y^m-1$, where $Q$
is a monic polynomial of degree $n$ with roots of multiplicities 
$p_1,\ldots,p_r$. In this case the curve $X=X_{\bar P}$ 
is the Riemann surface of the function $y=Q(x)^{1/m}$.
The number of components of this curve is the
greatest common divisor $d$ of $p_i$ and $m$. Also, the curve is a regular
covering of the line without $r$ points of degree $m$. Therefore,
the Euler characteristic of $X$ is $m(1-r)$, and thus 
$\dim B_2(A)=\dim H^1(X)=m(r-1)+d$. 

\end{example}

Let $P$ be a generic nonhomogeneous noncommutative
polynomial of degree $d$ in $n\ge 1$ variables. Let $A=A_n/(P)$.

\begin{proposition}
$\dim([A,A]\cap A[[A,A],A]A)/[[A,A],A]$ is $(d-1)^n$ if
$n$ is even, and $0$ if $n$ is odd.
\end{proposition}

\begin{proof}
Let $\bar P$ be the abelianization of $P$, and $X$
be the hypersurface defined by the equation $\bar P=0$
in $\C^n$. Then by Theorem \ref{pseudosm} and the results of \cite{EG},
the space $([A,A]\cap A[[A,A],A]A)/[[A,A],A]$ is isomorphic
to the odd cohomology $H^{\rm odd}(X)$.

Since $X$ is generic, it is obtained by removing of a smooth projective hypersurface of degree $d$ and dimension $n-2$ from one of degree $d$ and dimension $n-1$. Therefore, by the Lefschetz hyperplane section theorem,
$X$ has cohomology only in degrees $0$ and $n-1$. This implies the result in the case of odd $n$. If $n$ is even, the dimension of the odd cohomology is
$1-\chi$, where $\chi$ is the Euler characteristic of $X$.
So it remains to find $\chi$.

The computation of $\chi$ is well known, but we give it for the reader's convenience. We may assume that $X$ is the hypersurface
$X(d,n)$ defined by the equation 
$$
x_1^d+\cdots+x_n^d=1.
$$
Then by forgetting $x_n$ we get a degree $d$ surjective map
$X(d,n)\to \C^{n-1}$ which branches along $X(d,n-1)$ (where
there is $1$ instead of $d$ preimages). Thus if $\chi(d,n)$ denotes
the Euler characteristic of $X(d,n)$, then we have
$$
\chi(d,n)=d-(d-1)\chi(d,n-1).
$$
Since $\chi(d,1)=d$, we get by induction
$\chi(d,n)=1-(1-d)^n$. Hence the dimension in question is
$(d-1)^n$, as desired.
\end{proof}

{\bf Acknowledgments.} The author is very grateful
to B. Shoikhet for numerous very useful discussions, and in
particular for pointing out the relevance of the reference
\cite{CQ}. The work of the author was  partially supported by the NSF grant
 DMS-0504847.


\end{document}